\documentclass[12pt,russian,a4paper]{article}

\usepackage{amsmath,amssymb,amsthm,xspace,amscd}
\usepackage{amsfonts,amsxtra,latexsym}
\usepackage[cp1251]{inputenc}
\usepackage[russian]{babel}
\usepackage{geometry}

\usepackage{comment}

\usepackage{geometry}
\newtheorem{lemma}{Лемма}[section]
\newtheorem{theorem}{Теорема}[section]
\newtheorem{proposition}{Утверждение}[section]

\newtheorem{definition}{Определение}[section]

\begin{document}


\begin{center}
{\Large \bf Рост обобщенных алгебр Темперли-Либа} \\
13.11.2008
\end{center}

\section*{Введение}
$\;$

Алгебры Темперли-Либа и их $*\text{-}$представления изучались в ряде работ (\cite{tl}, \cite{jones} и др.) в связи
с моделями статистической физики. В \cite{graham} были введены обобщенные алгебры Темперли-Либа и среди них выделены конечномерные алгебры, связанные с графами $A_n,$ $D_n,$ $E_n$ и др.

В статье изучается рост обобщенных алгебр типа Темперли-Либа $TL_{\Gamma,\tau}.$  Изучаются их
размерности или, если алгебра бесконечномерная, рост.

\section{Класс алгебр $TL_\Gamma$}

\begin{definition}

Пусть $(W,S)$ является системой Кокстера с матрицей Кокстера $M=(m_{ij})_{i,j}$ и $q\in\mathbb C.$ Тогда алгеброй
Гекке $H_{W,S}$ будем называется унитальная ассоциативная алгебра над $\mathbb C$ порожденная образующими $t_s,$
$s\in S$ и соотношениями: $t_s^2=(q-1)t_s+q$ где $s\in S,$ $q\in\mathbb C$ и $t_rt_st_r\dots=t_st_rt_s\dots$ где в
каждой части равенства содержится $m_{ij}$ множителей, $r,s\in S,$ $r\neq s.$
\end{definition}

$H_{W,S}=\mathbb C\Bigl\langle t_1,\dots,t_n\mid t_k^2=(q-1)t_k+q,\;
(t_it_j)^{[\frac{m_{ij}}{2}]}t_i^{\sigma_{m_{ij}}}=(t_jt_i)^{[\frac{m_{ij}}{2}]}t_j^{\sigma_{m_{ij}}}\Bigr\rangle$

\begin{definition}
Обобщенной алгеброй Темперли-Либа $TL_{W,S},$ ассоциированной с системой Кокстера $(W,S)$ с графом $\Gamma$ будем
называть фактор-алгебру алгебры Гекке $H_{W,S},$ которая профакторизована по идеалу $I(\Gamma)$ порожденному
элементами $(t_it_j)^{[\frac{m_{ij}}{2}]}t_i^{\sigma_{m_{ij}}},$ где $(i,j)$ ребро в $\Gamma.$
\end{definition}

В дальнейшем обобщенную алгебру Темперли-Либа $TL_{W,S},$ ассоциированную с системой Кокстера $(W,S)$ с графом
$\Gamma$ будем обозначать через $TL_{\Gamma,\tau},$ где $\tau\in\mathbb C$ соответствующий параметр или
$TL_\Gamma$ и говорить, что она имеет тип $\Gamma.$ Также всюду ниже будем считать, что $\Gamma$ --- связный
неориентированный граф.

\begin{proposition}
Размерность и рост алгебры $TL_{\Gamma,\tau}$ не зависят от $\tau.$
\end{proposition}
\begin{proof}
Пусть $\Gamma_n$ --- связный неориентированный граф, $V\Gamma_n=n$. Достаточно показать, что старшие слова базиса
Гребнера алгебры $TL_{\Gamma_n}$ не зависят от $\tau.$ Доказательство проведем по индукции.

Напомним, что базисом Гребнера идеала $I$ называется множество элементов этого идеала $G\subset I$ такое, что для
всякого $g\in I$ старшее слово $g$ содержит в качестве подслова одно из старших слов элементов множества $G.$
Базис Гребнера $G$ называется минимальным, если никакое его собственное подмножество не является базисом Гребнера.
Алгоритм построения минимального базиса Гребнера в идеале, порожденном конечным множеством элементов свободной
алгебры, состоит из трех шагов: нормировка, редукция и композиция.

Предположим, что для алгебры $TL_{\Gamma_n}$ утверждение верно (основание индукции приведено ниже в примерах).
Соединим ребром граф $\Gamma_n$ с новой $n+1\text{-}$й вершиной (будем считать, что $(n,n+1)\in E\Gamma_n$).
Положим $\Gamma_{n+1}=\Gamma_n\cup\{n+1\}.$ Тогда алгебра $TL_{\Gamma_{n+1}}$ будет отличаться от $TL_{\Gamma_n}$
соотношениями $(t_nt_{n+1})^{[\frac{m_{ij}}{2}]}t_n^{\sigma_{m_{ij}}}=0,$
$(t_{n+1}t_n)^{[\frac{m_{ij}}{2}]}t_{n+1}^{\sigma_{m_{ij}}}=0$ и $t_it_{n+1}=t_{n+1}t_i$ для всех $i=1,\dots,n-1.$
Следовательно, $(t_nt_{n+1})^{[\frac{m_{ij}}{2}]}t_n^{\sigma_{m_{ij}}},$
$(t_{n+1}t_n)^{[\frac{m_{ij}}{2}]}t_{n+1}^{\sigma_{m_{ij}}}$ и $t_it_{n+1}-t_{n+1}t_i$ ($i=1,\dots,n-1$)
содержатся в базисе Гребнера алгебры $TL_{\Gamma_{n+1}},$ т.е старшие слова в этих базисных элементах не зависят
от $\tau.$ Но тогда по индуктивному предположению старшие слова новых элементов базиса Гребнера алгебры
$TL_{\Gamma_{n+1}}$ не будут зависеть от $\tau.$

\end{proof}


\section{Конечномерные алгебры $TL_{\Gamma}$} \label{dim<inf}
$\;$

\subsection{$TL_{A_n,\tau}$ ($n\geq 2$)}
Рассмотрим алгебру $TL_{A_n},$ ассоциированную с графом Дынкина $A_n$ ($n\geq 2$)

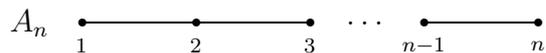
\setlength{\unitlength}{1mm}
\begin{picture}(80,20)(-10,20)
\linethickness{1pt} \thinlines \put(20,30){\line(1,0){30}} \put(65,30){\line(1,0){15}} \put(20,30){\circle*{1}}
\put(35,30){\circle*{1}} \put(50,30){\circle*{1}} \put(65,30){\circle*{1}} \put(80,30){\circle*{1}}
\put(13,30){\makebox(0,0)[a]{$A_n$}} \put(57.5,30){\makebox(0,0)[a]{$\dots$}} \put(20,27){\makebox(0,0)[a]{$_1$}}
\put(35,27){\makebox(0,0)[a]{$_2$}} \put(50,27){\makebox(0,0)[a]{$_3$}} \put(65,27){\makebox(0,0)[a]{$_{n-1}$}}
\put(80,27){\makebox(0,0)[a]{$_n$}}
\end{picture}

Алгебры $TL_{A_n,\tau}$ конечномерные как фактор-алгебры конечномерных алгебр Гекке, $\dim TL_{A_n}=\cfrac{1}{n+2}\begin{pmatrix}
  2n+2  \\
  n+1
\end{pmatrix}$ (см., например, \cite{jones}).

\subsection{$TL_{D_n}$ ($n\geq 4$)}

Рассмотрим граф Дынкина $D_n$ ($n\geq 4$)

\vspace{0.5cm}

\setlength{\unitlength}{1mm}
\begin{picture}(80,20)(-10,20)
\linethickness{1pt} \thinlines \put(35,30){\line(1,0){15}} \put(65,30){\line(1,0){15}}
\put(35,30){\line(-1,1){10}} \put(35,30){\line(-1,-1){10}} \put(25,20){\circle*{1}} \put(25,40){\circle*{1}}
\put(35,30){\circle*{1}} \put(50,30){\circle*{1}} \put(65,30){\circle*{1}} \put(80,30){\circle*{1}}
\put(18,30){\makebox(0,0)[a]{$D_n$}} \put(57.5,30){\makebox(0,0)[a]{$\dots$}} \put(22,20){\makebox(0,0)[a]{$_1$}}
\put(22,40){\makebox(0,0)[a]{$_2$}} \put(35,27){\makebox(0,0)[a]{$_3$}} \put(50,27){\makebox(0,0)[a]{$_4$}}
\put(65,27){\makebox(0,0)[a]{$_n$}} \put(80,27){\makebox(0,0)[a]{$_{n+1}$}}
\end{picture}

Алгебра $TL_{D_n,\tau}$ конечномерная, так как при $\tau=\frac 14$ является фактор-алгеброй конечномерной алгебры Гекке, ассоциированной с графом $D_n$ (см. например, \cite{burbaki}). Так $\dim TL_{D_4}=48,$ $\dim
TL_{D_5}=167,$ $\dim TL_{D_6}=593,$ $\dim TL_{D_7}=2144.$

\subsection{$TL_{E_n}$ ($n\geq 6)$}
Рассмотрим граф $E_n$ ($n\geq 6$)

\vspace{0.7cm} \setlength{\unitlength}{1mm}
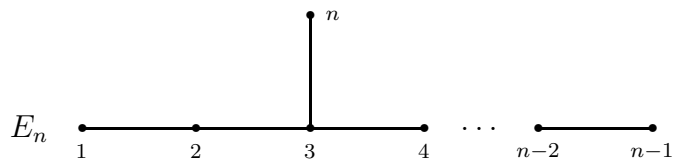
\begin{picture}(80,20)(-10,20)
\linethickness{1pt} \thinlines \put(20,30){\line(1,0){45}} \put(80,30){\line(1,0){15}} \put(50,30){\line(0,1){15}}
\put(20,30){\circle*{1}} \put(35,30){\circle*{1}} \put(50,30){\circle*{1}} \put(65,30){\circle*{1}}
\put(80,30){\circle*{1}} \put(95,30){\circle*{1}} \put(50,45){\circle*{1}} \put(13,30){\makebox(0,0)[a]{$E_n$}}
\put(72.5,30){\makebox(0,0)[a]{$\dots$}} \put(20,27){\makebox(0,0)[a]{$_1$}} \put(35,27){\makebox(0,0)[a]{$_2$}}
\put(50,27){\makebox(0,0)[a]{$_3$}} \put(65,27){\makebox(0,0)[a]{$_4$}} \put(80,27){\makebox(0,0)[a]{$_{n-2}$}}
\put(95,27){\makebox(0,0)[a]{$_{n-1}$}} \put(53,45){\makebox(0,0)[a]{$_n$}}

\end{picture}

Алгебры $TL_{E_n,\tau}$ при $n=6,7,8$ конечномерные, так как при $\tau=\frac 14$ являются фактор-алгебрами соответствующих конечномерных алгебр Гекке (см. например, \cite{burbaki}). Так $\dim TL_{E_6}=662,$ $\dim
TL_{E_7}=2670$ и $\dim TL_{E_8}=10846.$ Граф $E_9$ совпадает с расширенным графом Дынкина $\widetilde E_8.$ Алгебры $TL_{E_n,\tau}$ при $n\geq 9$ также конечномерные (см. \cite{graham}). Так $\dim TL_{E_9,\tau}=43409.$

Положим  $m_{ij}=4$ и $\sigma_{m_{ij}}=0$ при $(i,j)\in E\Gamma.$ Таким образом, мы будем рассматривать обобщенные
алгебры Темперли-Либа с соотношениями
$$\begin{cases}
p_ip_jp_ip_j=\tau p_ip_j, &\text{если}\; (i,j)_4\in E\Gamma;\\
p_jp_ip_jp_i=\tau p_jp_i, &\text{если}\; (i,j)_4\in E\Gamma;\\
p_ip_jp_i=\tau_i p_i, &\text{если}\; (i,j)\in E\Gamma;\\
p_jp_ip_j=\tau_i p_j, &\text{если}\; (i,j)\in E\Gamma;\\
p_ip_j=p_jp_i, &\text{если}\; (i,j)\not\in E\Gamma.
  \end{cases}
$$

\subsection{$TL_{B_n}$ ($n\geq 2$)}
Рассмотрим алгебру $TL_{B_n},$ ассоциированную с графом $B_n$

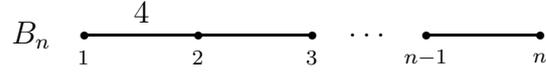
\setlength{\unitlength}{1mm}
\begin{picture}(80,20)(-10,20)
\linethickness{1pt} \thinlines \put(20,30){\line(1,0){30}} \put(65,30){\line(1,0){15}} \put(20,30){\circle*{1}}
\put(35,30){\circle*{1}} \put(50,30){\circle*{1}} \put(65,30){\circle*{1}} \put(80,30){\circle*{1}}
\put(13,30){\makebox(0,0)[a]{$B_n$}} \put(57.5,30){\makebox(0,0)[a]{$\dots$}} \put(27.5,33){\makebox(0,0)[a]{$4$}}
\put(20,27){\makebox(0,0)[a]{$_1$}} \put(35,27){\makebox(0,0)[a]{$_2$}} \put(50,27){\makebox(0,0)[a]{$_3$}}
\put(65,27){\makebox(0,0)[a]{$_{n-1}$}} \put(80,27){\makebox(0,0)[a]{$_n$}}

\end{picture}

Известно, что алгебры типа $B$ конечномерные (\cite{burbaki}, \cite{graham}). Например, $\dim TL_{B_2}=7,$ $\dim
TL_{B_3}=24,$ $\dim TL_{B_4}=83,$ $\dim TL_{B_5}=293.$

\subsection{$TL_{F_n}$ ($n\geq 4$)}

Рассмотрим алгебру $TL_{F_n},$ ассоциированную с графом $F_n$

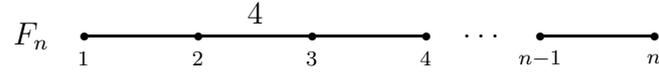
\setlength{\unitlength}{1mm}
\begin{picture}(80,20)(-10,20)
\linethickness{1pt} \thinlines \put(20,30){\line(1,0){45}} \put(80,30){\line(1,0){15}} \put(20,30){\circle*{1}}
\put(35,30){\circle*{1}} \put(50,30){\circle*{1}} \put(65,30){\circle*{1}} \put(80,30){\circle*{1}}
\put(95,30){\circle*{1}} \put(13,30){\makebox(0,0)[a]{$F_n$}} \put(72.5,30){\makebox(0,0)[a]{$\dots$}}
\put(42.5,33){\makebox(0,0)[a]{$4$}} \put(20,27){\makebox(0,0)[a]{$_1$}} \put(35,27){\makebox(0,0)[a]{$_2$}}
\put(50,27){\makebox(0,0)[a]{$_3$}} \put(65,27){\makebox(0,0)[a]{$_4$}} \put(80,27){\makebox(0,0)[a]{$_{n-1}$}}
\put(95,27){\makebox(0,0)[a]{$_n$}}

\end{picture}

Известно, что алгебры типа $F$ конечномерные (\cite{burbaki}, \cite{graham}). Например, $\dim TL_{F_4}=106,$ $\dim
TL_{F_5}=464,$ $\dim TL_{F_6}=2003.$

\subsection{$TL_{H_n}$ ($n\geq 2$)}

Пусть теперь $0\leq\tau_1,\tau_2\leq 1,$ $m_{ij}=4$ и $\sigma_{m_{ij}}=1$
при $(i,j)\in E\Gamma.$ Таким образом, мы будем рассматривать обобщенные алгебры Темперли-Либа с соотношениями
$$\begin{cases}
p_ip_jp_ip_jp_i=\tau_1p_ip_jp_i-\tau_2p_i, &\text{если}\; (i,j)_5\in E\Gamma;\\
p_jp_ip_jp_ip_j=\tau_1p_jp_ip_j-\tau_2p_j, &\text{если}\; (i,j)_5\in E\Gamma;\\
p_ip_jp_i=\tau_i p_i, &\text{если}\; (i,j)\in E\Gamma;\\
p_jp_ip_j=\tau_i p_j, &\text{если}\; (i,j)\in E\Gamma;\\
p_ip_j=p_jp_i, &\text{если}\; (i,j)\not\in E\Gamma.
  \end{cases}
$$

Рассмотрим алгебру $TL_{H_n},$ ассоциированную с графом $H_n$

\setlength{\unitlength}{1mm}
\begin{picture}(80,20)(-10,20)
\linethickness{1pt} \thinlines \put(20,30){\line(1,0){30}} \put(65,30){\line(1,0){15}} \put(20,30){\circle*{1}}
\put(35,30){\circle*{1}} \put(50,30){\circle*{1}} \put(65,30){\circle*{1}} \put(80,30){\circle*{1}}
\put(13,30){\makebox(0,0)[a]{$H_n$}} \put(57.5,30){\makebox(0,0)[a]{$\dots$}} \put(27.5,33){\makebox(0,0)[a]{$5$}}
\put(20,27){\makebox(0,0)[a]{$_1$}} \put(35,27){\makebox(0,0)[a]{$_2$}} \put(50,27){\makebox(0,0)[a]{$_3$}}
\put(65,27){\makebox(0,0)[a]{$_{n-1}$}} \put(80,27){\makebox(0,0)[a]{$_n$}}
\end{picture}

Известно (\cite{burbaki}, \cite{graham}), что алгебры типа $H$ конечномерные. Например, $\dim TL_{H_2}=9,$
$\dim TL_{H_3}=44,$ $\dim TL_{H_4}=195,$ $\dim TL_{H_5}=804.$

\subsection{$TL_{\widetilde G_2}$ и $TL_{l_2(p)}$ ($p=5,$ $p\geq 7$)}
Рассмотрим алгебру $TL_{\Gamma,\tau},$ ассоциированную с графом $\Gamma,$ где $\Gamma$ --- один из графов
$l_2(p)$ при $p=5,$ $p\geq7$ или $\widetilde G_2$

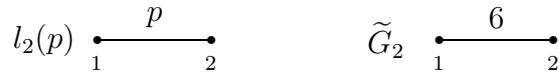
\setlength{\unitlength}{1mm}
\begin{picture}(80,20)(-10,20)
\linethickness{1pt} \thinlines \put(20,30){\line(1,0){15}} \put(65,30){\line(1,0){15}} \put(20,30){\circle*{1}}
\put(35,30){\circle*{1}} \put(65,30){\circle*{1}} \put(80,30){\circle*{1}} \put(13,30){\makebox(0,0)[a]{$l_2(p)$}}
\put(27.5,33){\makebox(0,0)[a]{$p$}} \put(58,30){\makebox(0,0)[a]{$\widetilde G_2$}}
\put(72.5,33){\makebox(0,0)[a]{$6$}} \put(20,27){\makebox(0,0)[a]{$_1$}} \put(35,27){\makebox(0,0)[a]{$_2$}}
\put(65,27){\makebox(0,0)[a]{$_1$}} \put(80,27){\makebox(0,0)[a]{$_2$}}
\end{picture}

Элементарные вычисления дают $\dim TL_{\widetilde G_2}=11$ и $\dim TL_{l_2(p)}=2p-1$ при $p=5,$
$p\geq7.$

\section{Алгебры $TL_\Gamma$ полиномиального роста}
Пусть $\mathcal A$ --- конечно порожденная ассоциативная алгебра и $\mathcal K(\mathcal A)$ --- семейство его
конечномерных подпространств. Говорят, что на $\mathcal A$ задана калибровка, если для любого $V\in\mathcal
K(\mathcal A)$ определена последовательность вложенных конечномерных подпространств
$$V^{(1)}\subseteq V^{(2)}\subseteq V^{(3)}\subseteq\dots, \; V^{(n)}\in\mathcal K(\mathcal A),$$ удовлетворяющая
следующему условию: если $V^{(k)}\subseteq W^{(m)},$ то для любого $n\in\mathbb N$ $V^{(kn)}\subseteq W^{(mn)}.$ С
каждым $V\in\mathcal K(\mathcal A)$ связана функция роста $d_V(n)=\dim V^{(n)}.$

На множестве таких функций определяются отношения предпорядка и эквивалентности: $f\leq g$ тогда и только тогда,
когда существуют такое $m\in\mathbb N$ и $c>0,$ что $f(n)\leq g(mn)$ для всех $n\in\mathbb N;$ $f\sim g$ тогда и
только тогда, когда $f\leq g$ и $g\leq f.$ Класс эквивалентности $f$ называется ростом $f$ и обоначается $[f].$

Для каждого конечномерного пространства $V$ алгебры $\mathcal A$ положим $V^{(n)}=V+V^1+\dots+V^n.$ Такая
калибровка определяет рост алгебры $\mathcal A.$

Все многочлены одинаковой степени $d$ имеют одинаковый рост $[n^d],$ который называется полиномиальным степени
$d.$ Если $[f]\leq [n^d]$ для некоторого $d,$ то рост $f$ считается полиномиальным. Рост $[2^n]$ называется
экспоненциальным.

Понятие роста определяется для любой монотонной функции, для алгебр и групп оказывается инвриантным понятием.

\subsection{$TL_{\widetilde A_n,\tau}$ ($n\geq 3$)} \label{221} Рассмотрим теперь расширенный граф Дынкина $\widetilde A_n$ ($n\geq 3$)

\vspace{0.3cm} \setlength{\unitlength}{1mm}
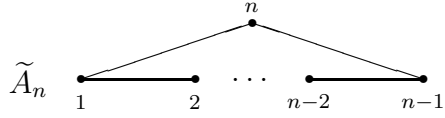
\begin{picture}(80,20)(-10,20)
\linethickness{1pt} \thinlines \put(35,30){\line(1,0){15}} \put(65,30){\line(1,0){15}}
\put(35,30){\line(3,1){22.5}} \put(80,30){\line(-3,1){22.5}} \put(35,30){\circle*{1}} \put(50,30){\circle*{1}}
\put(65,30){\circle*{1}} \put(80,30){\circle*{1}} \put(57.5,37.5){\circle*{1}}
\put(28,30){\makebox(0,0)[a]{$\widetilde A_n$}} \put(57.5,30){\makebox(0,0)[a]{$\dots$}}
\put(35,27){\makebox(0,0)[a]{$_1$}} \put(50,27){\makebox(0,0)[a]{$_2$}} \put(65,27){\makebox(0,0)[a]{$_{n-2}$}}
\put(80,27){\makebox(0,0)[a]{$_{n-1}$}} \put(57.5,39.5){\makebox(0,0)[a]{$_n$}}
\end{picture}

Известно (\cite{burbaki}), что группавая алгебра Кокстера, связанная с графом $\widetilde A_n$ бесконечномерна и имеет полиномиальный рост. Покажем, что алгебра и 
$TL_{\widetilde A_n}$ бесконечномерна и имеет линейный рост.

Рассмотрим алгебру $TL_{\widetilde A_n},$ ассоциированную с графом $\widetilde A_n$ ($n\geq 3$). Так как
максимальная длина старшего слова элементов базиса Гребнера алгебры $TL_{\widetilde A_n}$ равна $n,$ то вершинами
графа роста алгебры $TL_{\widetilde A_n}$ будут нормальные слова длины $n-1.$ Обозначим через $S_{n-1}$ множество
старших слов базиса Гребнера алгебры $TL_{\widetilde A_n}$ длины $n-1,$ т.е.
$S_{n-1}=\bigcup_{i_1,\dots,i_{n-1}}\prod_{k=1}^{n-1}p_{i_k}.$ Тогда граф роста будет иметь $n^{n-1}-|S_{n-1}|$
вершин.

Предположим, что алгебра $TL_{\widetilde A_n}$ имеет рост отличный от линейного. Тогда её граф роста будет
содержать не менее двух циклов соединенных путем. Обозначим через $l_1,$ $l_2$ два минимальных цикла, которые
соединеены путем. Пусть $t_1$ --- путь от $l_1$ к $l_2,$ а $t_2$ --- путь от $l_2$ к $l_1.$ Обозначим через $l_3$
минимальный цикл, соединяющий циклы $l_1$ и $l_2,$ который содержит пути $t_1$ и $t_2.$ Тогда для любых
$q_1,q_2\in l_3$ найдутся образующие $p_{i_1}$ и $p_{i_2}$ такие, что $q_1p_{i_1}=p_{i_2}q_2.$

Базис Гребнера алгебры $TL_{\widetilde A_n}$ содержит соотношения $p_1p_np_1-\tau p_1$ и $p_np_1p_n-\tau p_n,$
тогда $|Vl_i|=n$ при $i=1,2,3.$ Но граф роста алгебры $TL_{\widetilde A_n}$ не может иметь более двух циклов длины
$n.$ Тогда существуют такие вершины $q_1,q_2,q_3\in l_3$ и образующие $p_{i_1},$ $p_{i_2},$ $p_{i_3},$ $p_{i_4}$
для которых выполняются равенства $q_1p_{i_1}=p_{i_2}q_2$ и $p_{i_4}q_2=q_3p_{i_3},$ но тогда $TL_{\widetilde
A_n}$ имеет линейный рост.

Например, базис Гребнера алгебры $TL_{\widetilde A_3}$ состоит из элементов $\{ p_1^2-p_1,$ $p_2^2-p_2,$
$p_3^2-p_3,$ $p_1p_2p_1-\tau p_1,$ $p_2p_1p_2-\tau p_2,$ $p_2p_3p_2-\tau p_2,$ $p_3p_2p_3-\tau p_3,$
$p_1p_3p_1-\tau p_1,$ $p_3p_1p_3-\tau p_3\},$ а граф роста имеет вид

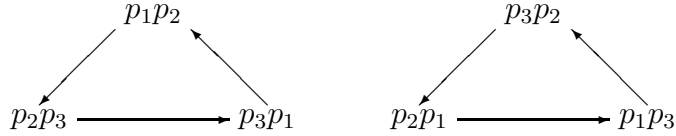
\setlength{\unitlength}{1mm}
\begin{picture}(80,20)(-10,20)
\linethickness{1pt} \thinlines \put(35,35){\makebox(0,0)[a]{$p_1p_2$}} \put(30,33){\vector(-1,-1){10}}
\put(20,21){\makebox(0,0)[a]{$p_2p_3$}} \put(25,21){\vector(1,0){20}} \put(50,21){\makebox(0,0)[a]{$p_3p_1$}}
\put(50,23){\vector(-1,1){10}} \put(85,35){\makebox(0,0)[a]{$p_3p_2$}} \put(80,33){\vector(-1,-1){10}}
\put(70,21){\makebox(0,0)[a]{$p_2p_1$}} \put(75,21){\vector(1,0){20}} \put(100,21){\makebox(0,0)[a]{$p_1p_3$}}
\put(100,23){\vector(-1,1){10}}
\end{picture}
\vspace{0.3cm}

\subsection{$TL_{\widetilde D_n}$ $(n\geq 4)$} \label{222}

Рассмотрим алгебру $TL_{\widetilde D_n},$ ассоциированную с графом $\widetilde D_n$ ($n\geq 4$)

\vspace{0.5cm} \setlength{\unitlength}{1mm}
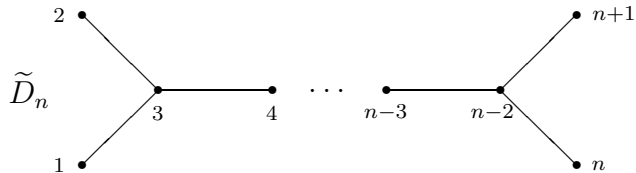
\begin{picture}(80,20)(-10,20)
\linethickness{1pt} \thinlines \put(35,30){\line(1,0){15}} \put(65,30){\line(1,0){15}}
\put(35,30){\line(-1,1){10}} \put(35,30){\line(-1,-1){10}} \put(80,30){\line(1,1){10}}
\put(80,30){\line(1,-1){10}} \put(35,30){\circle*{1}} \put(50,30){\circle*{1}} \put(65,30){\circle*{1}}
\put(80,30){\circle*{1}} \put(25,40){\circle*{1}} \put(25,20){\circle*{1}} \put(90,40){\circle*{1}}
\put(90,20){\circle*{1}} \put(18,30){\makebox(0,0)[a]{$\widetilde D_n$}} \put(57.5,30){\makebox(0,0)[a]{$\dots$}}
\put(22,20){\makebox(0,0)[a]{$_1$}} \put(22,40){\makebox(0,0)[a]{$_2$}} \put(35,27){\makebox(0,0)[a]{$_3$}}
\put(50,27){\makebox(0,0)[a]{$_4$}} \put(65,27){\makebox(0,0)[a]{$_{n-3}$}}
\put(79,27){\makebox(0,0)[a]{$_{n-2}$}} \put(93,20){\makebox(0,0)[a]{$_n$}}
\put(95,40){\makebox(0,0)[a]{$_{n+1}$}}
\end{picture}
\vspace{0.3cm}

Известно (\cite{burbaki}), что группавая алгебра Кокстера, связанная с графом $\widetilde D_n$ бесконечномерна и имеет полиномиальный рост. Аналогично пункту \ref{221} можно показать, что и алгебра $TL_{\widetilde D_n,\tau}$ бесклнечномерна и имеет линейный рост.

\subsection{$TL_{\widetilde E_6}$ и $TL_{\widetilde E_7}$} \label{223}
Рассмотрим граф Дынкина $\widetilde E_6$

\vspace{2.2cm}

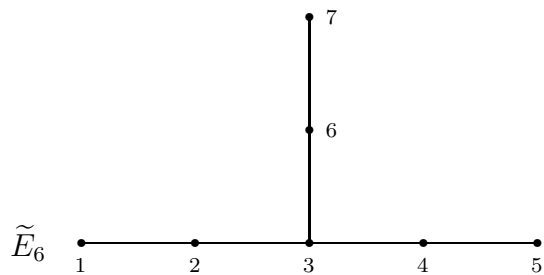
\setlength{\unitlength}{1mm}
\begin{picture}(80,20)(-10,20)
\linethickness{1pt} \thinlines \put(20,30){\line(1,0){60}} \put(50,30){\line(0,1){30}}

\put(20,30){\circle*{1}} \put(35,30){\circle*{1}} \put(50,30){\circle*{1}} \put(65,30){\circle*{1}}
\put(80,30){\circle*{1}}
\put(50,45){\circle*{1}} \put(50,60){\circle*{1}}

\put(13,30){\makebox(0,0)[a]{$\widetilde E_6$}}

\put(20,27){\makebox(0,0)[a]{$_1$}} \put(35,27){\makebox(0,0)[a]{$_2$}} \put(50,27){\makebox(0,0)[a]{$_3$}}
\put(65,27){\makebox(0,0)[a]{$_4$}} \put(80,27){\makebox(0,0)[a]{$_5$}} \put(53,45){\makebox(0,0)[a]{$_6$}}
\put(53,60){\makebox(0,0)[a]{$_7$}}
\end{picture}

Известно (\cite{burbaki}), что группавая алгебра Кокстера, связанная с графами $\widetilde E_6$ и $\widetilde E_8$ бесконечномерны и имеет полиномиальный рост. Аналогично пункту \ref{221} можно показать, что алгебры $TL_{\widetilde E_6,\tau}$ и $TL_{\widetilde E_7,\tau}$ бесклнечномерны и имеют линейный рост.ый рост.

\vspace{0.7cm} \setlength{\unitlength}{1mm}
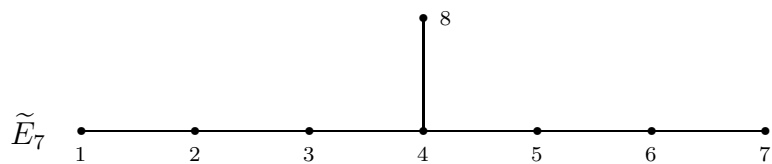
\begin{picture}(80,20)(-10,20)
\linethickness{1pt} \thinlines \put(20,30){\line(1,0){90}} \put(65,30){\line(0,1){15}} \put(20,30){\circle*{1}}
\put(35,30){\circle*{1}} \put(50,30){\circle*{1}} \put(65,30){\circle*{1}} \put(80,30){\circle*{1}}
\put(95,30){\circle*{1}} \put(110,30){\circle*{1}} \put(65,45){\circle*{1}}
\put(13,30){\makebox(0,0)[a]{$\widetilde E_7$}} \put(20,27){\makebox(0,0)[a]{$_1$}}
\put(35,27){\makebox(0,0)[a]{$_2$}} \put(50,27){\makebox(0,0)[a]{$_3$}} \put(65,27){\makebox(0,0)[a]{$_4$}}
\put(80,27){\makebox(0,0)[a]{$_5$}} \put(95,27){\makebox(0,0)[a]{$_6$}} \put(110,27){\makebox(0,0)[a]{$_7$}}
\put(68,45){\makebox(0,0)[a]{$_8$}}
\end{picture}

\subsection{$TL_{\widetilde A_1}$}
Рассмотрим алгебру $TL_{\widetilde A_1},$ ассоциированную с графом $\widetilde A_1$

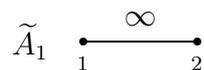
\setlength{\unitlength}{1mm}
\begin{picture}(80,20)(-10,20)
\linethickness{1pt} \thinlines \put(35,30){\line(1,0){15}} \put(35,30){\circle*{1}} \put(50,30){\circle*{1}}
\put(28,30){\makebox(0,0)[a]{$\widetilde A_1$}} \put(42.5,33){\makebox(0,0)[a]{$\infty$}}
\put(35,27){\makebox(0,0)[a]{$_1$}} \put(50,27){\makebox(0,0)[a]{$_2$}}
\end{picture}

Алгебра $TL_{\widetilde A_1}$ задается соотношениями
$$TL_{\widetilde A_1}=\mathbb C\Bigl\langle p_1,p_2\mid p_k^2=p_k^*=p_k,\; k=1,2\Bigr\rangle.$$

Так как число элементов линейного базиса алгебры $TL_{\widetilde A_1}$ фиксированной длины большей единицы
равно двум, то она имеет линейный рост (\cite{burbaki}).

\subsection{$TL_{\widetilde F_n}$ ($n\geq 6$)}
Рассмотрим алгебру $TL_{\widetilde F_n},$ ассоциированную с графом $\widetilde F_n$

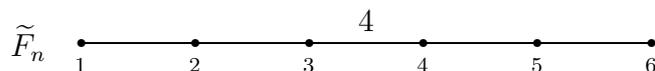
\setlength{\unitlength}{1mm}
\begin{picture}(80,20)(-10,20)
\linethickness{1pt} \thinlines \put(20,30){\line(1,0){75}} \put(20,30){\circle*{1}} \put(35,30){\circle*{1}}
\put(50,30){\circle*{1}} \put(65,30){\circle*{1}} \put(80,30){\circle*{1}} \put(95,30){\circle*{1}}
\put(20,27){\makebox(0,0)[a]{$_1$}} \put(35,27){\makebox(0,0)[a]{$_2$}} \put(50,27){\makebox(0,0)[a]{$_3$}}
\put(65,27){\makebox(0,0)[a]{$_4$}} \put(80,27){\makebox(0,0)[a]{$_5$}} \put(95,27){\makebox(0,0)[a]{$_6$}}
\put(13,30){\makebox(0,0)[a]{$\widetilde F_n$}} \put(57.5,33){\makebox(0,0)[a]{$4$}}
\end{picture}

Имеет линейный рост.

\subsection{$TL_{\widetilde B_n}$ ($n\geq 4$)}
 Рассмотрим алгебру $TL_{\widetilde B_n},$ ассоциированную с графом $\widetilde B_n$

\vspace{0.5cm}
\setlength{\unitlength}{1mm}
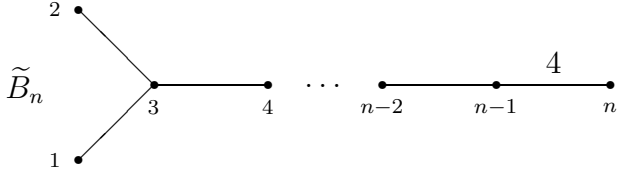
\begin{picture}(80,20)(-10,20)
\linethickness{1pt} \thinlines \put(20,30){\line(1,0){15}} \put(20,30){\line(-1,1){10}}
\put(20,30){\line(-1,-1){10}} \put(50,30){\line(1,0){30}} \put(20,30){\circle*{1}} \put(35,30){\circle*{1}}
\put(10,40){\circle*{1}} \put(10,20){\circle*{1}} \put(50,30){\circle*{1}} \put(65,30){\circle*{1}}
\put(80,30){\circle*{1}} \put(7,20){\makebox(0,0)[a]{$_1$}} \put(20,27){\makebox(0,0)[a]{$_3$}}
\put(7,40){\makebox(0,0)[a]{$_2$}} \put(35,27){\makebox(0,0)[a]{$_4$}} \put(50,27){\makebox(0,0)[a]{$_{n-2}$}}
\put(65,27){\makebox(0,0)[a]{$_{n-1}$}} \put(80,27){\makebox(0,0)[a]{$_n$}}
\put(3,30){\makebox(0,0)[a]{$\widetilde B_n$}} \put(42.5,30){\makebox(0,0)[a]{$\dots$}}
\put(72.5,33){\makebox(0,0)[a]{$4$}}
\end{picture}

Имеет линейный рост.

\subsection{$TL_{l_3(s)}$ ($s\geq 6$)}
Рассмотрим алгебру $TL_{l_3(p)},$ ассоциированную с графом $l_3(p)$ ($p\geq 6$)

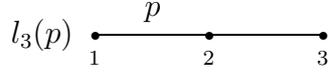
\setlength{\unitlength}{1mm}
\begin{picture}(80,20)(-10,20)
\linethickness{1pt} \thinlines \put(20,30){\line(1,0){30}} \put(20,30){\circle*{1}} \put(35,30){\circle*{1}}
\put(50,30){\circle*{1}} \put(13,30){\makebox(0,0)[a]{$l_3(p)$}} \put(27.5,33){\makebox(0,0)[a]{$p$}}
\put(20,27){\makebox(0,0)[a]{$_1$}} \put(35,27){\makebox(0,0)[a]{$_2$}} \put(50,27){\makebox(0,0)[a]{$_3$}}
\end{picture}

Имеет линейный рост.

\subsection{$TL_{\widetilde C_2}$}
Рассмотрим алгебру $TL_{\widetilde C_2},$ ассоциированную с графом $\widetilde C_2$

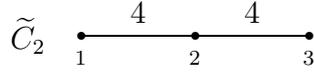
\setlength{\unitlength}{1mm}
\begin{picture}(80,20)(-10,20)
\linethickness{1pt} \thinlines \put(35,30){\line(1,0){30}} \put(35,30){\circle*{1}} \put(50,30){\circle*{1}}
\put(65,30){\circle*{1}} \put(28,30){\makebox(0,0)[a]{$\widetilde C_2$}} \put(35,27){\makebox(0,0)[a]{$_1$}}
\put(50,27){\makebox(0,0)[a]{$_2$}} \put(65,27){\makebox(0,0)[a]{$_3$}} \put(42.5,33){\makebox(0,0)[a]{$4$}}
\put(57.5,33){\makebox(0,0)[a]{$4$}}
\end{picture}

Известно (\cite{burbaki}), что алгебра $TL_{\widetilde C_2}$ имеет полиномиальный рост. Граф роста для
алгебры $TL_{\widetilde C_2}$ не содержит более одного цикла, соединенного путем, а значит его рост линейный.
Граф роста имеет вид

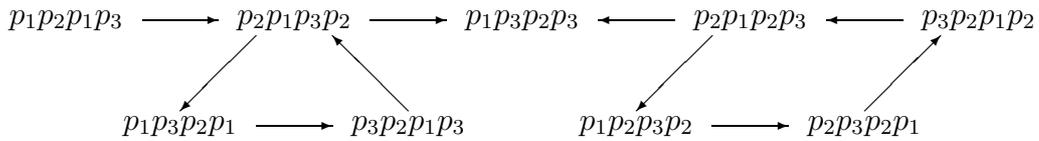
\setlength{\unitlength}{1mm}
\begin{picture}(80,20)(-10,20)
\linethickness{1pt} \thinlines \put(0,35){\makebox(0,0)[a]{$p_1p_2p_1p_3$}}
\put(30,35){\makebox(0,0)[a]{$p_2p_1p_3p_2$}} \put(25,33){\vector(-1,-1){10}} \put(45,23){\vector(-1,1){10}}
\put(60,35){\makebox(0,0)[a]{$p_1p_3p_2p_3$}} \put(90,35){\makebox(0,0)[a]{$p_2p_1p_2p_3$}}
\put(120,35){\makebox(0,0)[a]{$p_3p_2p_1p_2$}} \put(15,21){\makebox(0,0)[a]{$p_1p_3p_2p_1$}}
\put(25,21){\vector(1,0){10}} \put(45,21){\makebox(0,0)[a]{$p_3p_2p_1p_3$}}
\put(75,21){\makebox(0,0)[a]{$p_1p_2p_3p_2$}} \put(105,21){\makebox(0,0)[a]{$p_2p_3p_2p_1$}}
\put(85,21){\vector(1,0){10}} \put(10,35){\vector(1,0){10}} \put(40,35){\vector(1,0){10}}
\put(80,35){\vector(-1,0){10}} \put(110,35){\vector(-1,0){10}} \put(85,33){\vector(-1,-1){10}}
\put(105,23){\vector(1,1){10}}
\end{picture}
\vspace{0.3cm}

Базис Гребнера алгебры $TL_{\widetilde C_2}$ состоит из элементов $\{ p_1^2-p_1,$ $p_2^2-p_2,$ $p_3^2-p_3,$
$p_1p_2p_1p_2-\tau p_1p_2,$ $p_2p_1p_2p_1-\tau p_2p_1,$ $p_2p_3p_2p_3-\tau p_2p_3,$ $p_3p_2p_3p_2-\tau p_3p_2,$
$p_3p_1-p_1p_3,$ $p_2p_3p_2p_1p_3-\tau p_2p_1p_3,$ $p_1p_3p_2p_1p_3-\tau p_1p_3p_2 \}.$

\subsection{$TL_{\widetilde C_n}$ ($n\geq 3$)}
2. Рассмотрим алгебру $TL_{\widetilde C_n},$ ассоциированную с графом $\widetilde C_n$ ($n\geq 3$)


\setlength{\unitlength}{1mm}
\begin{picture}(80,20)(-10,20)
\linethickness{1pt} \thinlines \put(20,30){\line(1,0){30}} \put(65,30){\line(1,0){30}} \put(20,30){\circle*{1}}
\put(35,30){\circle*{1}} \put(50,30){\circle*{1}} \put(65,30){\circle*{1}} \put(80,30){\circle*{1}}
\put(95,30){\circle*{1}} \put(13,30){\makebox(0,0)[a]{$\widetilde C_n$}} \put(57.5,30){\makebox(0,0)[a]{$\dots$}}
\put(20,27){\makebox(0,0)[a]{$_1$}} \put(35,27){\makebox(0,0)[a]{$_2$}} \put(50,27){\makebox(0,0)[a]{$_3$}}
\put(65,27){\makebox(0,0)[a]{$_{n-2}$}} \put(80,27){\makebox(0,0)[a]{$_{n-1}$}}
\put(95,27){\makebox(0,0)[a]{$_n$}} \put(27.5,33){\makebox(0,0)[a]{$4$}} \put(87.5,33){\makebox(0,0)[a]{$4$}}
\end{picture}

Известно (\cite{burbaki}), что алгебры типа $\widetilde C_n$ имеют линейный рост.

\section{Алгебры $TL_\Gamma$ экспоненциального роста}

\subsection{}

Рассмотрим алгебру $TL_{K_6},$ ассоциированную с графом $K_6$

\vspace{0.7cm} \setlength{\unitlength}{1mm}
\begin{picture}(80,20)(-10,20)
\linethickness{1pt} \thinlines \put(35,30){\line(1,0){30}} \put(50,30){\line(0,1){15}} \put(50,30){\line(1,1){10}}
\put(50,30){\line(-1,1){10}} \put(35,30){\circle*{1}} \put(50,30){\circle*{1}} \put(65,30){\circle*{1}}
\put(50,45){\circle*{1}} \put(60,40){\circle*{1}} \put(40,40){\circle*{1}} \put(28,30){\makebox(0,0)[a]{$K_6$}}
\put(35,27){\makebox(0,0)[a]{$_1$}} \put(50,27){\makebox(0,0)[a]{$_2$}} \put(65,27){\makebox(0,0)[a]{$_3$}}
\put(63,40){\makebox(0,0)[a]{$_4$}} \put(37,40){\makebox(0,0)[a]{$_6$}} \put(53,45){\makebox(0,0)[a]{$_5$}}
\end{picture}

Обозначим через $F=\mathbb C\Bigl\langle x_1,x_2\Bigr\rangle$ свободную алгебру, порожденную двумя образующими.
Положим $q_1=p_2p_3p_4p_2p_1p_5$ и $q_2=p_2p_3p_4p_2p_1p_6.$ Линейный базис $TL_{K_6}$ состоит из всевозможных
слов, в которых не содержатся следующие подслова $p_1^2,$ $p_2^2,$ $p_3^2,$ $p_4^2,$ $p_5^2,$ $p_6^2,$
$p_1p_2p_1,$ $p_2p_1p_2,$ $p_2p_3p_2,$ $p_3p_2p_3,$ $p_2p_4p_2,$ $p_4p_2p_4,$ $p_2p_5p_2,$ $p_5p_2p_5,$
$p_2p_6p_2,$ $p_6p_2p_6,$ $p_3p_1,$ $p_4p_1,$ $p_5p_1,$ $p_6p_1,$ $p_4p_3,$ $p_5p_3,$ $p_6p_3,$ $p_5p_4,$
$p_6p_4,$ $p_6p_5,$ $p_3p_2p_1p_3,$ $p_4p_2p_1p_4,$ $p_4p_2p_3p_4,$ $p_5p_2p_1p_5,$ $p_5p_2p_3p_5,$
$p_5p_2p_4p_5,$ $p_6p_2p_1p_6,$ $p_6p_2p_3p_6,$ $p_6p_2p_4p_6,$ $p_6p_2p_5p_6,$ $p_1p_3p_2p_1,$ $p_1p_4p_2p_1,$
$p_1p_5p_2p_1,$ $p_1p_6p_2p_1,$ $p_3p_4p_2p_3,$ $p_3p_5p_2p_3,$ $p_3p_6p_2p_3,$ $p_4p_5p_2p_4,$ $p_4p_6p_2p_4,$
$p_5p_6p_2p_5,$ $p_1p_3p_4p_2p_1,$ $p_1p_3p_5p_2p_1,$ $p_1p_3p_6p_2p_1,$ $p_1p_4p_5p_2p_1,$ $p_1p_4p_6p_2p_1,$
$p_1p_5p_6p_2p_1,$ $p_3p_4p_2p_1p_3,$ $p_3p_4p_5p_2p_3,$ $p_3p_4p_6p_2p_3,$ $p_3p_5p_2p_1p_3,$ $p_3p_5p_6p_2p_3,$
$p_3p_6p_2p_1p_3,$ $p_4p_5p_2p_1p_4,$ $p_4p_5p_2p_3p_4,$ $p_4p_5p_6p_2p_4,$ $p_4p_6p_2p_1p_4,$ $p_4p_6p_2p_3p_4,$
$p_5p_6p_2p_1p_5,$ $p_5p_6p_2p_3p_5,$ $p_5p_6p_2p_4p_5,$ $p_4p_2p_1p_3p_4,$ $p_5p_2p_1p_3p_5,$ $p_5p_2p_1p_4p_5,$
$p_5p_2p_3p_4p_5,$ $p_6p_2p_1p_3p_6,$ $p_6p_2p_1p_4p_6,$ $p_6p_2p_1p_5p_6,$ $p_6p_2p_3p_4p_6,$ $p_6p_2p_3p_5p_6,$
$p_6p_2p_4p_5p_6,$ $p_1p_3p_4p_5p_2p_1,$ $p_1p_3p_4p_6p_2p_1,$ $p_1p_3p_5p_6p_2p_1,$ $p_1p_4p_5p_6p_2p_1,$
$p_3p_4p_5p_2p_1p_3,$ $p_3p_4p_5p_6p_2p_3,$ $p_3p_4p_6p_2p_1p_3,$ $p_3p_5p_6p_2p_1p_3,$ $p_4p_5p_6p_2p_1p_4,$
$p_4p_5p_6p_2p_3p_4,$ $p_4p_5p_2p_1p_3p_4,$ $p_4p_6p_2p_1p_3p_4,$ $p_5p_6p_2p_1p_3p_5,$ $p_5p_6p_2p_1p_4p_5,$
$p_5p_6p_2p_3p_4p_5,$ $p_5p_2p_1p_3p_4p_5,$ $p_6p_2p_1p_3p_4p_6,$ $p_6p_2p_1p_3p_5p_6,$ $p_6p_2p_1p_4p_5p_6,$
$p_6p_2p_3p_4p_5p_6,$ $p_1p_3p_4p_5p_6p_2p_1,$ $p_3p_4p_5p_6p_2p_1p_3,$ $p_4p_5p_6p_2p_1p_3p_4,$
$p_5p_6p_2p_1p_3p_4p_5,$ $p_6p_2p_1p_3p_4p_5p_6.$

Пусть $\varphi : F \longrightarrow G$ гомоморфизм $\varphi(x_i)=q_i,$ где $G=\mathbb C\Bigl\langle
q_1,q_2\Bigr\rangle$ подалгебра алгебры $TL_{K_6,\tau}.$ Так как ни одно из запрещенных слов не содержится во
всевозможных комбинациях элементов $q_1$ и $q_2,$ то $\varphi$ является изоморфизмом. Следовательно, в алгебре
$TL_{K_6,\tau}$ существует свободная подалгебра $G$ порожденная двумя образующими, тогда $TL_{K_6,\tau}$ имеет
экспоненциальный рост.

\subsection{}
Рассмотрим граф

\vspace{0.5cm}

\setlength{\unitlength}{1mm}
\begin{picture}(80,20)(-10,20)
\linethickness{1pt} \thinlines \put(35,30){\line(1,0){45}} \put(50,30){\line(-1,1){10}}
\put(50,30){\line(-1,-1){10}} \put(35,30){\circle*{1}} \put(50,30){\circle*{1}} \put(65,30){\circle*{1}}
\put(80,30){\circle*{1}} \put(40,40){\circle*{1}} \put(40,20){\circle*{1}} \put(35,27){\makebox(0,0)[a]{$_1$}}
\put(50,27){\makebox(0,0)[a]{$_2$}} \put(65,27){\makebox(0,0)[a]{$_3$}} \put(80,27){\makebox(0,0)[a]{$_4$}}
\put(37,20){\makebox(0,0)[a]{$_5$}} \put(37,40){\makebox(0,0)[a]{$_6$}}
\end{picture}

Соответствующая алгебра $TL_\Gamma$ имеет экспоненциальный рост, так как у неё существует свободная подалгебра
порожденная двумя образующими (например, можно взять подалгебру порожденную элементами
$q_1=p_2p_1p_3p_2p_4p_3p_5p_2p_1p_6p_2p_3p_4p_5p_2p_1p_3p_2p_5p_6$ и
$q_2=p_2p_1p_3p_2p_4p_3p_6p_2p_1p_5p_2p_3p_4p_6p_2p_1p_3p_2p_5p_6$).

\subsection{}

Рассмотрим граф

\vspace{0.5cm} \setlength{\unitlength}{1mm}
\begin{picture}(80,20)(-10,20)
\linethickness{1pt} \thinlines \put(50,30){\line(1,0){15}} \put(50,30){\line(-1,1){10}}
\put(50,30){\line(-1,-1){10}} \put(65,30){\line(1,1){10}} \put(65,30){\line(1,-1){10}} \put(75,20){\line(1,0){15}}
\put(50,30){\circle*{1}} \put(65,30){\circle*{1}} \put(40,40){\circle*{1}} \put(40,20){\circle*{1}}
\put(75,40){\circle*{1}} \put(75,20){\circle*{1}} \put(90,20){\circle*{1}} \put(50,27){\makebox(0,0)[a]{$_3$}}
\put(65,27){\makebox(0,0)[a]{$_4$}} \put(37,40){\makebox(0,0)[a]{$_2$}} \put(37,20){\makebox(0,0)[a]{$_1$}}
\put(78,40){\makebox(0,0)[a]{$_5$}} \put(75,17){\makebox(0,0)[a]{$_6$}} \put(90,17){\makebox(0,0)[a]{$_7$}}
\end{picture}
\vspace{0.5cm}

Соответствующая алгебра $TL_\Gamma$ имеет экспоненциальный рост, так как у неё существует свободная подалгебра
порожденная двумя образующими (например, можно взять подалгебру порожденную элементами
$q_1=p_7p_6p_4p_3p_1p_2p_3p_4p_5p_6p_4p_3p_1p_2p_3p_4p_5$ и $q_2=p_6p_4p_3p_1p_2p_3p_4p_5$).

\subsection{}
Рассмотрим граф

\vspace{2.2cm} \setlength{\unitlength}{1mm}
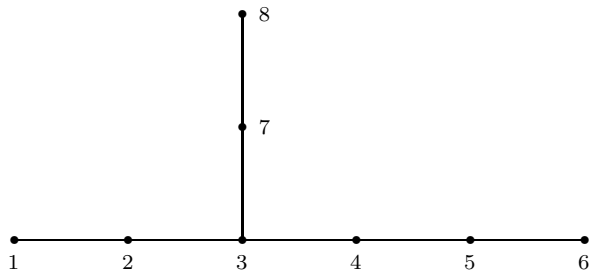
\begin{picture}(80,20)(-10,20)
\linethickness{1pt} \thinlines \put(20,30){\line(1,0){75}} \put(50,30){\line(0,1){30}} \put(20,30){\circle*{1}}
\put(35,30){\circle*{1}} \put(50,30){\circle*{1}} \put(65,30){\circle*{1}} \put(80,30){\circle*{1}}
\put(95,30){\circle*{1}} \put(50,45){\circle*{1}} \put(50,60){\circle*{1}} \put(20,27){\makebox(0,0)[a]{$_1$}}
\put(35,27){\makebox(0,0)[a]{$_2$}} \put(50,27){\makebox(0,0)[a]{$_3$}} \put(65,27){\makebox(0,0)[a]{$_4$}}
\put(80,27){\makebox(0,0)[a]{$_5$}} \put(95,27){\makebox(0,0)[a]{$_6$}} \put(53,45){\makebox(0,0)[a]{$_7$}}
\put(53,60){\makebox(0,0)[a]{$_8$}}
\end{picture}

Соответствующая алгебра $TL_\Gamma$ имеет экспоненциальный рост, так как у неё существует свободная подалгебра
порожденная двумя образующими (например, можно взять подалгебру порожденную элементами
$q_1=p_2p_7p_3p_4p_5p_8p_7p_3p_2p_1p_4p_3$ и \\
$q_2=p_2p_7p_3p_4p_5p_6p_8p_7p_3p_2p_4p_3p_5p_4p_7p_3p_2p_1p_8p_7p_3p_2p_4p_3p_5p_4p_6p_5p_7p_3p_2p_1p_4p_3$).

\subsection{}
Рассмотрим граф

\vspace{0.7cm} \setlength{\unitlength}{1mm}
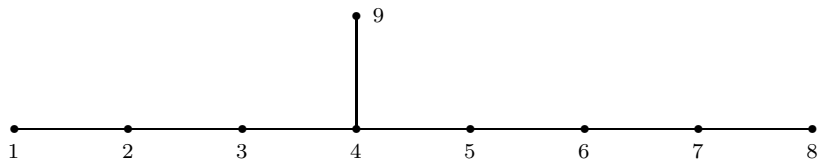
\begin{picture}(80,20)(-10,20)
\linethickness{1pt} \thinlines \put(20,30){\line(1,0){105}} \put(65,30){\line(0,1){15}} \put(20,30){\circle*{1}}
\put(35,30){\circle*{1}} \put(50,30){\circle*{1}} \put(65,30){\circle*{1}} \put(80,30){\circle*{1}}
\put(95,30){\circle*{1}} \put(110,30){\circle*{1}} \put(125,30){\circle*{1}} \put(65,45){\circle*{1}}
\put(20,27){\makebox(0,0)[a]{$_1$}} \put(35,27){\makebox(0,0)[a]{$_2$}} \put(50,27){\makebox(0,0)[a]{$_3$}}
\put(65,27){\makebox(0,0)[a]{$_4$}} \put(80,27){\makebox(0,0)[a]{$_5$}} \put(95,27){\makebox(0,0)[a]{$_6$}}
\put(110,27){\makebox(0,0)[a]{$_7$}} \put(125,27){\makebox(0,0)[a]{$_8$}} \put(68,45){\makebox(0,0)[a]{$_9$}}
\end{picture}

Соответствующая алгебра $TL_\Gamma$ имеет экспоненциальный рост, так как у неё существует свободная подалгебра
порожденная двумя образующими (например, можно взять подалгебру порожденную элементами
$q_1=p_4p_5p_6p_8p_4p_1p_2p_5p_7p_8p_4p_8p_9p_4p_2p_3$ и $q_2=r_1r_2,$ где
$r_1=p_4p_5p_6p_7p_8p_9p_5p_4p_3p_5p_4p_6p_7p_8p_4p_1p_3p_7p_8p_9p_4p_2$ и \\
$r_2=p_5p_7p_9p_1p_5p_4p_2p_8p_4p_6p_7p_4p_5p_7p_9p_4p_1p_3p_5p_7p_8p_4p_2p_3p_4p_5p_6p_8p_4p_1p_2p_3$).

\subsection{}
Рассмотрим граф

\vspace{0.5cm} \setlength{\unitlength}{1mm}
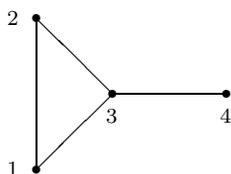
\begin{picture}(80,20)(-10,20)
\linethickness{1pt} \thinlines \put(50,30){\line(1,0){15}} \put(50,30){\line(-1,1){10}}
\put(50,30){\line(-1,-1){10}} \put(40,20){\line(0,1){20}} \put(50,30){\circle*{1}} \put(65,30){\circle*{1}}
\put(40,40){\circle*{1}} \put(40,20){\circle*{1}}
\put(50,27){\makebox(0,0)[a]{$_3$}} \put(65,27){\makebox(0,0)[a]{$_4$}} \put(37,20){\makebox(0,0)[a]{$_1$}}
\put(37,40){\makebox(0,0)[a]{$_2$}}
\end{picture}

Соответствующая алгебра $TL_\Gamma$ имеет экспоненциальный рост, так как у неё существует свободная подалгебра
порожденная двумя образующими (например, можно взять подалгебру порожденную элементами $q_1=p_2p_1p_3$ и
$q_2=p_2p_1p_3p_2p_1p_4p_3$).

\subsection{}
Рассмотрим алгебру $TL_\Gamma,$ ассоциированную с графом $\Gamma$

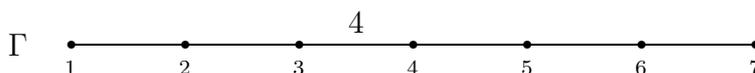
\setlength{\unitlength}{1mm}
\begin{picture}(80,20)(-10,20)
\linethickness{1pt} \thinlines \put(20,30){\line(1,0){90}} \put(20,30){\circle*{1}} \put(35,30){\circle*{1}}
\put(50,30){\circle*{1}} \put(65,30){\circle*{1}} \put(80,30){\circle*{1}} \put(95,30){\circle*{1}}
\put(110,30){\circle*{1}} \put(20,27){\makebox(0,0)[a]{$_1$}} \put(35,27){\makebox(0,0)[a]{$_2$}}
\put(50,27){\makebox(0,0)[a]{$_3$}} \put(65,27){\makebox(0,0)[a]{$_4$}} \put(80,27){\makebox(0,0)[a]{$_5$}}
\put(95,27){\makebox(0,0)[a]{$_6$}} \put(110,27){\makebox(0,0)[a]{$_7$}} \put(13,30){\makebox(0,0)[a]{$\Gamma$}}
\put(57.5,33){\makebox(0,0)[a]{$4$}}
\end{picture}

Алгебра $TL_\Gamma$ имеет экспоненциальный рост, так как у неё существует свободная подалгебра порожденная
двумя образующими (например, можно взять подалгебру порожденную элементами
$q_1=p_2p_3p_1p_6p_5p_4p_3p_2p_4p_3p_5p_4$ и \\ $q_2=p_3p_6p_5p_4p_3p_2p_1p_7p_6p_5p_4p_3p_2p_4p_3p_5p_4$).

\subsection{}
Рассмотрим алгебру $TL_\Gamma,$ ассоциированную с графом $\Gamma$

\vspace{0.5cm}
\setlength{\unitlength}{1mm}
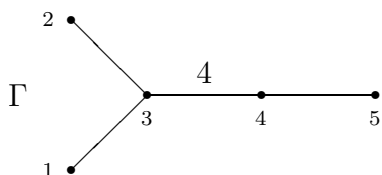
\begin{picture}(80,20)(-10,20)
\linethickness{1pt} \thinlines \put(35,30){\line(1,0){30}} \put(35,30){\line(-1,1){10}}
\put(35,30){\line(-1,-1){10}} \put(35,30){\circle*{1}} \put(50,30){\circle*{1}} \put(25,40){\circle*{1}}
\put(25,20){\circle*{1}} \put(65,30){\circle*{1}} \put(22,20){\makebox(0,0)[a]{$_1$}}
\put(35,27){\makebox(0,0)[a]{$_3$}} \put(22,40){\makebox(0,0)[a]{$_2$}} \put(50,27){\makebox(0,0)[a]{$_4$}}
\put(65,27){\makebox(0,0)[a]{$_5$}} \put(18,30){\makebox(0,0)[a]{$\Gamma$}} \put(42.5,33){\makebox(0,0)[a]{$4$}}
\end{picture}

Алгебра $TL_\Gamma$ имеет экспоненциальный рост, так как у неё существует свободная подалгебра порожденная
двумя образующими (например, можно взять подалгебру порожденную элементами $q_1=p_1p_2p_3p_4p_3$ и
$q_2=p_5p_4p_3p_1p_2p_3p_4p_3$).

\subsection{}
Рассмотрим алгебру

\vspace{0.5cm}
\setlength{\unitlength}{1mm}
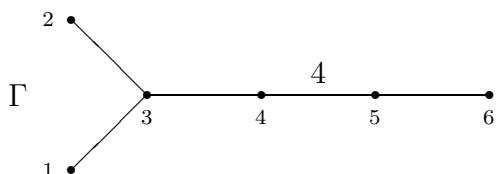
\begin{picture}(80,20)(-10,20)
\linethickness{1pt} \thinlines \put(35,30){\line(1,0){45}} \put(35,30){\line(-1,1){10}}
\put(35,30){\line(-1,-1){10}} \put(35,30){\circle*{1}} \put(50,30){\circle*{1}} \put(25,40){\circle*{1}}
\put(25,20){\circle*{1}} \put(65,30){\circle*{1}} \put(80,30){\circle*{1}} \put(22,20){\makebox(0,0)[a]{$_1$}}
\put(35,27){\makebox(0,0)[a]{$_3$}} \put(22,40){\makebox(0,0)[a]{$_2$}} \put(50,27){\makebox(0,0)[a]{$_4$}}
\put(65,27){\makebox(0,0)[a]{$_5$}} \put(80,27){\makebox(0,0)[a]{$_6$}} \put(18,30){\makebox(0,0)[a]{$\Gamma$}}
\put(57.5,33){\makebox(0,0)[a]{$4$}}
\end{picture}

(например, можно взять подалгебру порожденную элементами $q_1=p_2p_6p_5p_4p_3p_1p_5p_4p_3$ и
$q_2=p_2p_5p_4p_3p_1p_5p_4p_3$).

\subsection{}
Рассмотрим алгебру

\vspace{0.5cm}
\setlength{\unitlength}{1mm}
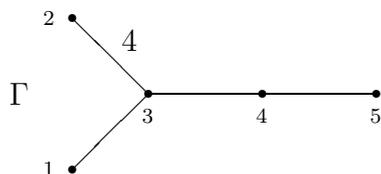
\begin{picture}(80,20)(-10,20)
\linethickness{1pt} \thinlines \put(35,30){\line(1,0){30}} \put(35,30){\line(-1,1){10}}
\put(35,30){\line(-1,-1){10}} \put(35,30){\circle*{1}} \put(50,30){\circle*{1}} \put(25,40){\circle*{1}}
\put(25,20){\circle*{1}} \put(65,30){\circle*{1}} \put(22,20){\makebox(0,0)[a]{$_1$}}
\put(35,27){\makebox(0,0)[a]{$_3$}} \put(22,40){\makebox(0,0)[a]{$_2$}} \put(50,27){\makebox(0,0)[a]{$_4$}}
\put(65,27){\makebox(0,0)[a]{$_5$}} \put(18,30){\makebox(0,0)[a]{$\Gamma$}} \put(32.5,37){\makebox(0,0)[a]{$4$}}
\end{picture}

(например, можно взять подалгебру порожденную элементами $q_1=p_2p_5p_4p_3p_1p_2p_3p_2p_4p_3$ и
$q_2=p_1p_2p_3p_2p_4p_3$).

\subsection{}
Рассмотрим алгебру

\vspace{1.2cm} \setlength{\unitlength}{1mm}
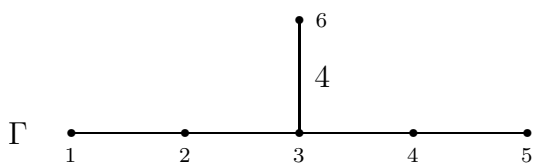
\begin{picture}(80,20)(-10,20)
\linethickness{1pt} \thinlines \put(20,30){\line(1,0){60}} \put(50,30){\line(0,1){15}} \put(20,30){\circle*{1}}
\put(35,30){\circle*{1}} \put(50,30){\circle*{1}} \put(65,30){\circle*{1}} \put(80,30){\circle*{1}}
\put(50,45){\circle*{1}} \put(13,30){\makebox(0,0)[a]{$\Gamma$}} \put(20,27){\makebox(0,0)[a]{$_1$}}
\put(35,27){\makebox(0,0)[a]{$_2$}} \put(50,27){\makebox(0,0)[a]{$_3$}} \put(65,27){\makebox(0,0)[a]{$_4$}}
\put(80,27){\makebox(0,0)[a]{$_5$}} \put(53,45){\makebox(0,0)[a]{$_6$}} \put(53,37.5){\makebox(0,0)[a]{$4$}}

\end{picture}

(например, можно взять подалгебру порожденную элементами $q_1=p_2p_1p_4p_3p_2p_5p_4p_3p_5p_4p_5p_6p_3$ и
$q_2=p_2p_4p_3p_5p_4p_5p_6p_3$).

\subsection{}
Рассмотрим алгебру

\vspace{1.2cm} \setlength{\unitlength}{1mm}
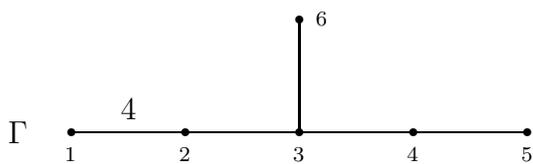
\begin{picture}(80,20)(-10,20)
\linethickness{1pt} \thinlines \put(20,30){\line(1,0){60}} \put(50,30){\line(0,1){15}} \put(20,30){\circle*{1}}
\put(35,30){\circle*{1}} \put(50,30){\circle*{1}} \put(65,30){\circle*{1}} \put(80,30){\circle*{1}}
\put(50,45){\circle*{1}} \put(13,30){\makebox(0,0)[a]{$\Gamma$}} \put(20,27){\makebox(0,0)[a]{$_1$}}
\put(35,27){\makebox(0,0)[a]{$_2$}} \put(50,27){\makebox(0,0)[a]{$_3$}} \put(65,27){\makebox(0,0)[a]{$_4$}}
\put(80,27){\makebox(0,0)[a]{$_5$}} \put(53,45){\makebox(0,0)[a]{$_6$}} \put(27.5,33){\makebox(0,0)[a]{$4$}}

\end{picture}

(например, можно взять подалгебру порожденную элементами $q_1=p_3p_2p_1p_4p_3p_2p_1p_5p_4p_3p_2p_1p_6$ и
$q_2=p_3p_2p_1p_4p_3p_2p_1p_5p_6p_3p_2p_1p_4p_3p_2p_1p_6$).

\subsection{}
Рассмотрим алгебру

\vspace{1.2cm}

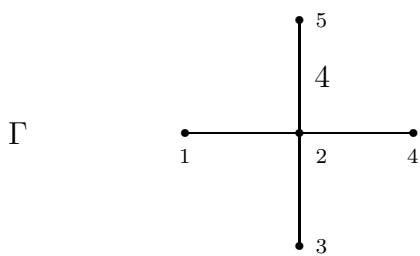
\setlength{\unitlength}{1mm}
\begin{picture}(80,20)(-10,20)
\linethickness{1pt} \thinlines \put(35,30){\line(1,0){30}} \put(50,15){\line(0,1){30}} \put(35,30){\circle*{1}}
\put(50,30){\circle*{1}} \put(65,30){\circle*{1}} \put(50,15){\circle*{1}} \put(50,45){\circle*{1}}
\put(13,30){\makebox(0,0)[a]{$\Gamma$}} \put(35,27){\makebox(0,0)[a]{$_1$}}
\put(53,27){\makebox(0,0)[a]{$_2$}} \put(65,27){\makebox(0,0)[a]{$_4$}} \put(53,15){\makebox(0,0)[a]{$_3$}}
\put(53,45){\makebox(0,0)[a]{$_5$}} \put(53,37.5){\makebox(0,0)[a]{$4$}}

\end{picture}

\vspace{1cm} (например, можно взять подалгебру порожденную элементами $q_1=p_2p_1p_3p_2p_4p_5$ и
$q_2=p_2p_3p_5p_2p_4p_5$).

\subsection{}
Рассмотрим алгебру

\vspace{1.2cm}

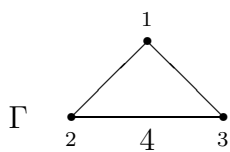
\setlength{\unitlength}{1mm}
\begin{picture}(80,20)(-10,20)
\linethickness{1pt} \thinlines \put(50,30){\line(1,0){20}} \put(50,30){\line(1,1){10}}
\put(70,30){\line(-1,1){10}} \put(50,30){\circle*{1}} \put(70,30){\circle*{1}} \put(60,40){\circle*{1}}
\put(43,30){\makebox(0,0)[a]{$\Gamma$}} \put(50,27){\makebox(0,0)[a]{$_2$}}
\put(70,27){\makebox(0,0)[a]{$_3$}} \put(60,43){\makebox(0,0)[a]{$_1$}} \put(60,27){\makebox(0,0)[a]{$4$}}
\end{picture}

(например, можно взять подалгебру порожденную элементами $q_1=p_1p_2p_3p_2p_1p_3p_2p_3$ и $q_2=p_1p_2p_3$).

\subsection{}
Рассмотрим алгебру

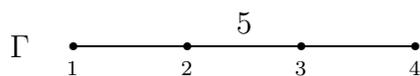
\setlength{\unitlength}{1mm}
\begin{picture}(80,20)(-10,20)
\linethickness{1pt} \thinlines \put(20,30){\line(1,0){45}} \put(20,30){\circle*{1}} \put(35,30){\circle*{1}}
\put(50,30){\circle*{1}} \put(65,30){\circle*{1}} \put(13,30){\makebox(0,0)[a]{$\Gamma$}}
\put(42.5,33){\makebox(0,0)[a]{$5$}} \put(20,27){\makebox(0,0)[a]{$_1$}} \put(35,27){\makebox(0,0)[a]{$_2$}}
\put(50,27){\makebox(0,0)[a]{$_3$}} \put(65,27){\makebox(0,0)[a]{$_4$}}
\end{picture}

(например, можно взять подалгебру порожденную элементами $q_1=p_2p_1p_3p_2p_3p_2p_4p_3$ и
$q_2=p_2p_1p_3p_2p_2p_4p_3$).

\subsection{}
Рассмотрим алгебру

\vspace{1.2cm}
\setlength{\unitlength}{1mm}
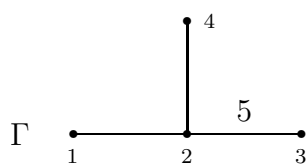
\begin{picture}(80,20)(-10,20)
\linethickness{1pt} \thinlines \put(20,30){\line(1,0){30}} \put(35,30){\line(0,1){15}} \put(20,30){\circle*{1}}
\put(35,30){\circle*{1}} \put(50,30){\circle*{1}} \put(35,45){\circle*{1}} \put(13,30){\makebox(0,0)[a]{$\Gamma$}}
\put(42.5,33){\makebox(0,0)[a]{$5$}} \put(20,27){\makebox(0,0)[a]{$_1$}} \put(35,27){\makebox(0,0)[a]{$_2$}}
\put(50,27){\makebox(0,0)[a]{$_3$}} \put(38,45){\makebox(0,0)[a]{$_4$}}
\end{picture}

(например, можно взять подалгебру порожденную элементами $q_1=p_2p_1p_3p_2p_3p_2p_1p_4p_2p_3$ и
$q_2=p_2p_1p_3p_2p_3p_4$).

\subsection{}
Рассмотрим алгебру

\vspace{1.2cm}

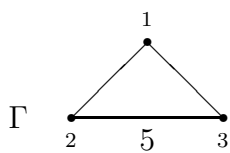
\setlength{\unitlength}{1mm}
\begin{picture}(80,20)(-10,20)
\linethickness{1pt} \thinlines \put(50,30){\line(1,0){20}} \put(50,30){\line(1,1){10}}
\put(70,30){\line(-1,1){10}} \put(50,30){\circle*{1}} \put(70,30){\circle*{1}} \put(60,40){\circle*{1}}
\put(43,30){\makebox(0,0)[a]{$\Gamma$}} \put(50,27){\makebox(0,0)[a]{$_2$}}
\put(70,27){\makebox(0,0)[a]{$_3$}} \put(60,43){\makebox(0,0)[a]{$_1$}} \put(60,27){\makebox(0,0)[a]{$5$}}
\end{picture}

(например, можно взять подалгебру порожденную элементами $q_1=p_1p_2p_3p_2p_1p_3p_2p_1p_3p_2p_3$ и
$q_2=p_1p_2p_3$).

\subsection{}
Рассмотрим алгебру

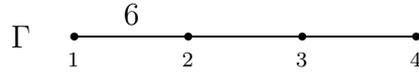
\setlength{\unitlength}{1mm}
\begin{picture}(80,20)(-10,20)
\linethickness{1pt} \thinlines \put(20,30){\line(1,0){45}} \put(20,30){\circle*{1}} \put(35,30){\circle*{1}}
\put(50,30){\circle*{1}} \put(65,30){\circle*{1}} \put(13,30){\makebox(0,0)[a]{$\Gamma$}}
\put(27.5,33){\makebox(0,0)[a]{$6$}} \put(20,27){\makebox(0,0)[a]{$_1$}} \put(35,27){\makebox(0,0)[a]{$_2$}}
\put(50,27){\makebox(0,0)[a]{$_3$}} \put(65,27){\makebox(0,0)[a]{$_4$}}
\end{picture}

(например, можно взять подалгебру порожденную элементами $q_1=p_1p_2p_1p_3p_2$ и
$q_2=p_1p_4p_3p_2p_1p_2p_1p_3p_2$).

\subsection{}
Рассмотрим алгебру

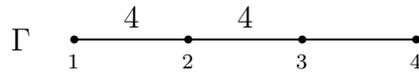
\setlength{\unitlength}{1mm}
\begin{picture}(80,20)(-10,20)
\linethickness{1pt} \thinlines \put(20,30){\line(1,0){45}} \put(20,30){\circle*{1}} \put(35,30){\circle*{1}}
\put(50,30){\circle*{1}} \put(65,30){\circle*{1}} \put(13,30){\makebox(0,0)[a]{$\Gamma$}}
\put(20,27){\makebox(0,0)[a]{$_1$}} \put(35,27){\makebox(0,0)[a]{$_2$}} \put(50,27){\makebox(0,0)[a]{$_3$}}
\put(65,27){\makebox(0,0)[a]{$_4$}} \put(27.5,33){\makebox(0,0)[a]{$4$}} \put(42.5,33){\makebox(0,0)[a]{$4$}}
\end{picture}

(например, можно взять подалгебру порожденную элементами $q_1=p_1p_3p_2$ и \\
$q_2=p_1p_3p_2p_4p_3p_2p_1p_2p_3p_2p_1p_4p_3p_2$).

\subsection{}
Рассмотрим алгебру

\vspace{1cm} \setlength{\unitlength}{1mm}
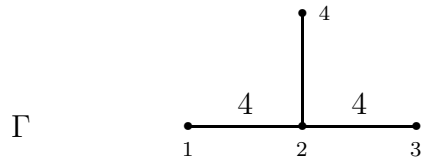
\begin{picture}(80,20)(-10,20)
\linethickness{1pt} \thinlines

\put(35,30){\line(1,0){30}} \put(50,30){\line(0,1){15}}

\put(35,30){\circle*{1}} \put(50,30){\circle*{1}} \put(65,30){\circle*{1}} \put(50,45){\circle*{1}}

\put(13,30){\makebox(0,0)[a]{$\Gamma$}}

\put(35,27){\makebox(0,0)[a]{$_1$}} \put(50,27){\makebox(0,0)[a]{$_2$}} \put(65,27){\makebox(0,0)[a]{$_3$}}
\put(53,45){\makebox(0,0)[a]{$_4$}}

\put(42.5,33){\makebox(0,0)[a]{$4$}} \put(57.5,33){\makebox(0,0)[a]{$4$}}
\end{picture}

(например, можно взять подалгебру порожденную элементами $q_1=p_1p_3p_2$ и $q_2=p_1p_3p_2p_1p_3p_4p_2$).

\subsection{}
Рассмотрим алгебру

\vspace{1cm} \setlength{\unitlength}{1mm}
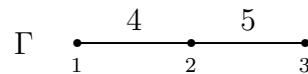
\begin{picture}(80,20)(-10,20)
\linethickness{1pt} \thinlines \put(35,30){\line(1,0){30}} \put(35,30){\circle*{1}} \put(50,30){\circle*{1}}
\put(65,30){\circle*{1}} \put(28,30){\makebox(0,0)[a]{$\Gamma$}} \put(35,27){\makebox(0,0)[a]{$_1$}}
\put(50,27){\makebox(0,0)[a]{$_2$}} \put(65,27){\makebox(0,0)[a]{$_3$}} \put(42.5,33){\makebox(0,0)[a]{$4$}}
\put(57.5,33){\makebox(0,0)[a]{$5$}}
\end{picture}

(например, можно взять подалгебру порожденную элементами $q_1=p_1p_3p_2$ и $q_2=p_1p_3p_2p_3p_2p_1p_2p_3p_2$).

\subsection{}

Рассмотрим алгебру

\vspace{1cm} \setlength{\unitlength}{1mm}
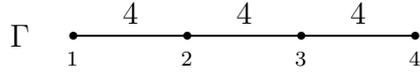
\begin{picture}(80,20)(-10,20)
\linethickness{1pt} \thinlines \put(20,30){\line(1,0){45}} \put(20,30){\circle*{1}} \put(35,30){\circle*{1}}
\put(50,30){\circle*{1}} \put(65,30){\circle*{1}} \put(13,30){\makebox(0,0)[a]{$\Gamma$}}
\put(20,27){\makebox(0,0)[a]{$_1$}} \put(35,27){\makebox(0,0)[a]{$_2$}} \put(50,27){\makebox(0,0)[a]{$_3$}}
\put(65,27){\makebox(0,0)[a]{$_4$}} \put(27.5,33){\makebox(0,0)[a]{$4$}} \put(42.5,33){\makebox(0,0)[a]{$4$}}
\put(57.5,33){\makebox(0,0)[a]{$4$}}
\end{picture}

(например, можно взять подалгебру порожденную элементами $q_1=p_1p_2p_1p_3p_4$ и $q_2=p_1p_2p_4p_2p_4$).

\section{Связь между ростом графа и ростом его подграфа.}
Существует связь между ростом графа и ростом его подграфа.

\begin{lemma}
Если к графу $\Gamma$ добавить новую вершину и соединить её ребром с одной из вершин $\Gamma,$ то размерность или
рост соответствующей алгебры $TL_\Gamma$ не уменьшается.
\end{lemma}
\begin{proof}
Пусть $\Gamma$ --- связный неориентированный граф. Добавив к $\Gamma$ вершину и соединив ребром с одной из вершин
$\Gamma,$ то получим новый граф $\Gamma'.$ Но тогда $TL_{\Gamma,\tau}$ является подалгеброй $TL_{\Gamma',\tau}$ и
$\dim TL_{\Gamma',\tau}\geq \dim TL_{\Gamma,\tau}.$
\end{proof}

\begin{lemma}
Если в графе $\Gamma$ соединить ребром любые две вершины, то размерность или рост соответствующей алгебры
$TL_\Gamma$ не уменьшается.
\end{lemma}
\begin{proof}
Так как мы рассматриваем только простые графы без петель, то достаточно будет рассмотреть граф $\widetilde A_3,$
полученный из $A_3$ добавлением нового ребра. Алгебра $TL_{A_3,\tau}$ конечномерная, а $TL_{\widetilde A_3,\tau}$
--- бесконечномерная.
\end{proof}

Пусть $x_1$ произвольная точка графа $\Gamma$ валентности $m.$ Разделим множество вершин $V\Gamma,$ которые
соединены с вершиной $x_1$ ребром на два помножества $V_1\Gamma$ и $V_2\Gamma.$ Построим новый граф $\Gamma'.$
Возьмем вершину $x_1$ и новую вершину $x_2$ (где $x_2\not\in V\Gamma$) и соединим их ребром. Все вершины из
$V_1\Gamma$ соединим с $x_1,$ а вершины из $V_2\Gamma$ соединим с $x_2.$ Остальные вешины соединим как в исходном
графе $\Gamma.$ Такую процедуру будем называть как сделать из точки ребро.

\begin{lemma}
Если в графе $\Gamma$ из точки сделать ребро, то размерность или рост соответствующей алгебры $TL_\Gamma$ не
уменьшается.
\end{lemma}
\begin{proof}
Для доказательства достаточно рассмотреть граф $\widetilde D_5,$ полученный из $\widetilde D_4$ прибавлением ребра
вместо точки. Алгебры $TL_{\widetilde D_4,\tau}$ и $TL_{\widetilde D_5,\tau}$ бесконечномерны и имеют линейный
рост.
\end{proof}

\section{Основная теорема}

\begin{theorem}
Пусть $\Gamma$ --- связный неориентированный граф, тогда
\begin{itemize}
\item[---] алгебра $TL_\Gamma$ конечномерная тогда и только тогда, когда граф $\Gamma$ совпадает с одним из графов $A_n,$ $D_n,$ $E_n,$ $B_n,$ $F_n,$ $H_n,$ $\widetilde G_2$ и $l_2(p)$ ($p=5,$ $p\geq 7$);
\item[---] алгебра $TL_\Gamma$ имеет линейный рост тогда и только тогда, когда граф $\Gamma$ совпадает с одним из графов $\widetilde A_n,$ $\widetilde D_n,$ $\widetilde E_6,$ $\widetilde E_7,$ $\widetilde A_1,$ $\widetilde F_n,$ $\widetilde B_n,$ $\widetilde C_n$ и $l_3(s)$ ($s\geq 6$);
\item[---] алгебра $TL_\Gamma$ имеет экспоненциальный рост тогда и только тогда, когда граф $\Gamma$ не совпадает ни с одним из графов $A_n,$ $D_n,$ $E_n,$ $\widetilde A_n,$ $\widetilde D_n,$ $\widetilde E_6,$ $\widetilde E_7,$ $B_n,$ $F_n,$ $H_n,$ $\widetilde G_2,$ $l_2(p)$ ($p=5,$ $p\geq 7$), $\widetilde A_1,$ $\widetilde F_n,$ $\widetilde B_n,$ $\widetilde C_n$ и $l_3(s)$ ($s\geq 6$).
\end{itemize}

\end{theorem}

\end{document}